\documentclass[11pt,leqno]{amsart}
\usepackage{amsmath}
\usepackage{amsthm}
\usepackage{amssymb}
\theoremstyle{definition}

\theoremstyle{plain}

\begin{document}
\setcounter{page}{171}
\begin{flushright}
\scriptsize{{\em Journal of the Nigerian Association of Mathematical Physics}
\\
{\em Volume} {\bf 10}, ({\em November} 2006), 171--176 (Original Print)}
\\
\copyright{\em J. of NAMP}
\end{flushright}
\vspace{16mm}

\begin{center}
{\normalsize\bf Some remarks on certain Bazilevic functions} \\[12mm]
     {\sc K. O. BABALOLA$^{1}$} \\ [8mm]

\begin{minipage}{123mm}
{\small {\sc Abstract.}
In this note we give some sufficient conditions for an analytic function $f(z)$ normalized by $f'(0)=1$ to belong to certain subfamilies of the class of Bazilevic functions. In earlier works, the closure property of many classes of functions under the Bernardi
integral have been considered. The converse of this problem is also considered here.}
\end{minipage}
\end{center}

 \renewcommand{\thefootnote}{}
\footnotetext{$^1$\,Department of Mathematics, University of Ilorin, Ilorin, Nigeria. abuqudduus\symbol{64}yahoo.com}

\def\iff{if and only if }
\def\S{Smarandache }
\newcommand{\norm}[1]{\left\Vert#1\right\Vert}

\vskip 12mm

{\bf 1.0 Introduction }
\medskip

Let $C$ be the complex plane. Denote by $A$ the class of functions:
$$f(z)=z+a_2z^2+\cdots\eqno(1.1)$$
which are analytic in the unit disk $E=\{z\colon|z|<1\}$.
Let $S$ be the subclass of $A$ consisting only of univalent functions in $E$.
Bazilevic \cite{4} isolated a subclass of $S$, which consists of functions defined
by the integral
$$f(z)=\left\{\frac{\alpha}{1+\xi^2}\int_0^z[p(\nu)-i\xi]\nu^{-\left(1+\frac{i\alpha\xi}{1+\xi^2}\right)}g(\nu)^{\left(\frac{\alpha}{1+\xi^2}\right)}d\nu\right\}^{\frac{1+i\xi}{\alpha}}\eqno(1.2)$$
where $p\in P$ (consisting of analytic functions $p(z)$ which have positive real part in
$E$ and normalized by $p(0)=1$) and $g\in S^\ast$ (the subclass of $S$ satisfying Re
$zg'(z)/g(z)>0$, i. e. starlike in $E$). The numbers $\alpha>0$ and $\xi$ are real and all powers
are meant as principal determinations only. The family of functions (1.2), which we
denote by $B(\alpha,\xi,p,g)$ is known as family of Bazilevic functions. He proved that
every Bazilevic function is univalent in $E$. Apart from this, very little is known
regarding the family as a whole. However, with some simplifications, it has been
possible to understand and investigate the family. Indeed it is easy to verify that
with special choices of the parameters $\alpha$ and $\xi$, and the function $g(z)$, the family
$B(\alpha,\xi,p,g)$ comes down to some well-known subclasses of univalent functions. For
instance, if we put $\xi=0$, we have
$$f(z)=\left\{\alpha\int_0^z\frac{p(\nu)}{\nu}g(\nu)^\alpha d\nu\right\}^{\frac{1}{\alpha}}.\eqno(1.3)$$

On differentiation, the expression (1.3) yields
$$\frac{zf'(z)f(z)^{\alpha-1}}{g(z)^\alpha}=p(z),\;\;z\in E\eqno(1.4)$$
that is,
$$Re\frac{zf'(z)f(z)^{\alpha-1}}{g(z)^\alpha}>0,\;\;z\in E.\eqno(1.5)$$

The subclasses of Bazilevic functions satisfying (1.5) are called Bazilevic
functions of type $\alpha$ and are denoted by $B(\alpha)$ (see \cite{8}). Noonan \cite{6} gave a plausible
description of functions of the class $B(\alpha)$ as those functions in $S$ for which for each
$r<1$, the tangent to the curve $U_\alpha(r)={f(re^{i\theta})^\alpha:0\leq\theta<2\pi}$ never turns back on
itself as much as $\pi$ radian.\vskip 2mm

With specific choices of the associated parameters the class $B(\alpha)$ reduces to
the well-known families of close-to-convex, starlike and convex functions.
Furthermore, if $g(z)=z$ in (1.5) we have the family $B_1(\alpha)$ \cite{8} of functions
satisfying:
$$Re\frac{zf'(z)f(z)^{\alpha-1}}{z^\alpha}>0,\;\;z\in E.\eqno(1.6)$$

Abdulhalim \cite{1} introduced and studied a generalization of functions
satisfying (1.6) as follows:
$$Re\frac{D^nf(z)^\alpha}{z^\alpha}>0,\;\;z\in E\eqno(1.7)$$
where the operator $D^n$ is the well-known Salagean derivative operator defined by
the relations $D^0f(z)=f(z)$ and $D^nf(z)=z[D^{n-1}f(z)]'$ \cite{1}. He denoted this class of
functions by $B_n(\alpha)$. It is immediately obvious that Abdulhalim's generalization has
extraneously included analytic functions satisfying:
$$Re\frac{f(z)^\alpha}{z^\alpha}>0,\;\;z\in E\eqno(1.8)$$
which are largely non-univalent in the unit disk. The case $\alpha=1$ here, coincides
with the class of functions studied by Yamaguchi \cite{9}.\vskip 2mm

Opoola \cite{7} further generalized functions defined by the geometric condition
(1.7) by choosing a real number $\beta$ $(0\leq\beta<1)$ such that:
$$Re\frac{D^nf(z)^\alpha}{\alpha^nz^\alpha}>\beta,\;\;z\in E\eqno(1.9)$$
(see \cite{2,7}).\vskip 2mm

In this article we provide sufficient conditions for an analytic function $f\in A$
to belong to certain subfamilies of the class of Bazilevic functions. In section 2 we
state some preliminary lemmas. We prove the main results in section 3.

 \medskip

{\bf 2.0 Preliminary Lemmas}\vskip 2mm

{\bf Definition 2.1.}\vskip 2mm

Let $u=u_1+u_2i$, $v=v_1+v_2i$. Define $\Psi$ as the set of functions $\psi(u,v):C\times C\to C$ satisfying:
 
{\rm(a)} $\psi(u,v)$ is continuous in a domain $\Omega$ of $C\times C$,

{\rm(b)} $(1,0)\in\Omega$ and Re$\psi(1,0)>0$,

{\rm(c)} Re$\psi(u_2i, v_1)\leq 0$ when $(u_2i, v_1)\in\Omega$ and $2v_1\leq -(1+u_2^2)$.\vskip 2mm

The following two examples of the set $\Psi$ add to those mentioned in \cite{2,5}.

{\rm(i)} $\psi_1(u,v)=v/\xi u$, $\xi>0$ is real and $\Omega=C\times C$.

{\rm(ii)} $\psi_2(u,v)=\tfrac{1}{2}+v/(\xi(1+u))$, $)<\xi\leq 1$ and $\Omega=[C-\{-1\}]\times C$.
\medskip

{\bf Definition 2.2.}\vskip 2mm

Let $\psi\in\Psi$ with corresponding domain $\Omega$. Define $P(\Psi)$ as the set of functions $p(z)$ given as $p(z)=1+p_1z+ p_2z^2+...$ which are regular in $E$ and satisfy:
 
{\rm(i)} $(p(z),zp^\prime(z))\in\Omega$

{\rm(ii)} Re$\psi(p(z),zp^\prime(z))>0$ when $z\in E$.

The above definitions are abridged forms of the more general concepts discussed in \cite{2}.\vskip 2mm

{\bf Lemma 2.3.}(\cite{2})\vskip 2mm

{\em Let $p\in P(\Psi)$. Then Re $p(z)>0$}.\vskip 2mm

{\bf Lemma 2.4.}(\cite{3})\vskip 2mm

{\em Let $f\in A$ and $\alpha>0$ be real. If $D^{n+1}f(z)^\alpha/D^nf(z)^\alpha$ is independent of $n$ for $z\in E$, then
\[
\frac{D^{n+1}f(z)^\alpha}{D^nf(z)^\alpha}=\alpha\frac{D^{n+1}f(z)}{D^nf(z)}.
\]
}

{\bf Lemma 2.5.}(\cite{8})\vskip 2mm

{\em Let $p\in P$. Then
$$|zp'(z)|\leq\frac{2r}{1-r^2}Re\;p(z),\;\;|z|=r.$$
The result is sharp. Equality holds for the function $p(z)=(1+z)/(1-z)$}.\vskip 2mm

\medskip

{\bf 3.0 Main Results}\vskip 2mm

{\bf Theorem 3.1.}\vskip 2mm

{\em Let $f\in A$. If $\alpha>0$ and Re$\left(\tfrac{D^{n+1}f(z)}{D^nf(z)}-1\right)>0$, then $f\in B_n(\alpha)$}.

{\bf Proof.}\vskip 2mm

Let $p(z)=\tfrac{D^nf(z)}{\alpha^nz^\alpha}$. Then we have $\tfrac{1}{\alpha}\tfrac{D^{n+1}f(z)^\alpha}{D^nf(z)^\alpha}-1=\tfrac{zp'(z)}{\alpha p(z)}$. Thus applying Lemma 2.4 we obtain $\tfrac{D^{n+1}f(z)}{D^nf(z)}-1=\tfrac{zp'(z)}{\alpha p(z)}=\psi(p(z),zp'(z))$, where $\psi=\psi_1(u,v)=\xi v/u$ (with $\xi=1/\alpha$) belongs to $\Psi$. Hence, by Lemma 2.3, we have the implication that Re $(p(z),zp(z))>0\Rightarrow$ Re $p(z)>0$, that is, Re$\tfrac{D^nf(z)}{z^\alpha}>0$, that is, $f(z)$ belongs to the class $B_n(\alpha)$.\vskip 2mm

{\bf Corollary 3.2.}\vskip 2mm

{\em Let $f\in A$. If $\alpha>0$, then

{\rm(a)} $$Re\left(\frac{zf'(z)}{f(z)}-1\right)>0\Rightarrow Re\frac{f(z)^\alpha}{z^\alpha}>0,$$

{\rm(b)} $$Re\frac{zf''(z)}{f'(z)}>0\Rightarrow Re\frac{f(z)^{\alpha-}f'(z)}{z^{\alpha-1}}>0.$$}

In particular if $\alpha=1$, we have\vskip 2mm

{\bf Corollary 3.3.}\vskip 2mm

{\em Let $f\in A$. Then

{\rm(a)} $$Re\left(\frac{zf'(z)}{f(z)}-1\right)>0\Rightarrow Re\frac{f(z)}{z}>0,$$

{\rm(b)} $$Re\frac{zf''(z)}{f'(z)}>0\Rightarrow Re f'(z)>0.$$}

{\bf Remark 3.4.}\vskip 2mm

In fact for $0\leq\beta<1$, $\psi_1(u,v)=v/\xi u$, with $\xi>0$ real and $\Omega=[C-{0}]\times C$ is
found to be contained in the set $\Psi_\beta$, $0\leq\beta<1$, defined in \cite{2}. Thus the above
theorem and corollaries can be extended easily to the class $T_n^\alpha(\beta)$.\vskip 2mm

{\bf Theorem 3.5.}\vskip 2mm

{\em Let $f\in A$. If $0<\alpha\leq 1$, then

$$Re\frac{D^{n+1}f(z)}{D^nf(z)}>\frac{1}{2}\Rightarrow Re\frac{D^nf(z)^\alpha}{\alpha^nz^\alpha}>\frac{1}{2}.$$}

{\bf Proof.}\vskip 2mm

Let $p(z)=2\tfrac{D^nf(z)}{\alpha^nz^\alpha}-1$. Then we have $\tfrac{1}{\alpha}\tfrac{D^{n+1}f(z)^\alpha}{D^nf(z)^\alpha}=1+\tfrac{zp'(z)}{\alpha(1+ p(z))}$. Thus applying Lemma 2.4 we have
$$\frac{D^{n+1}f(z)}{D^nf(z)}-\frac{1}{2}=\frac{1}{2}+\frac{zp'(z)}{\alpha p(z)}=\psi(p(z),zp'(z)),$$
where $\psi_2(u,v)=\tfrac{1}{2}+v/(\xi(1+u))$ (with $\xi=\alpha$) and  $\Omega=[C-{-1}]\times C$. Therefor, by Lemma 2.3, we have the implication that Re $(p(z),zp(z))>0\Rightarrow$ Re $p(z)>0$, that is,
$$Re\frac{D^{n+1}f(z)}{D^nf(z)}>\frac{1}{2}\Rightarrow Re\frac{D^nf(z)^\alpha}{\alpha^nz^\alpha}>\frac{1}{2}.$$

{\bf Corollary 3.6.}\vskip 2mm

{\em Let $f\in A$. If $0<\alpha\leq 1$, then

{\rm(a)} $$Re\frac{zf'(z)}{f(z)}>\frac{1}{2}\Rightarrow Re\frac{f(z)^\alpha}{z^\alpha}>\frac{1}{2},$$

{\rm(b)} $$Re\left\{1+\frac{zf''(z)}{f'(z)}\right\}>\frac{1}{2}\Rightarrow Re\frac{f(z)^{\alpha-}f'(z)}{z^{\alpha-1}}>\frac{1}{2}.$$}

Furthermore, if $\alpha=1$ we have\vskip 2mm

{\bf Corollary 3.7.}\vskip 2mm

{\em Let $f\in A$. Then

{\rm(a)} $$Re\frac{zf'(z)}{f(z)}>\frac{1}{2}\Rightarrow Re\frac{f(z)}{z}>\frac{1}{2},$$

{\rm(b)} $$Re\left\{1+\frac{zf''(z)}{f'(z)}\right\}>\frac{1}{2}\Rightarrow Re f'(z)>\frac{1}{2}.$$}

{\bf Remark 3.8.}\vskip 2mm

The first part of Corollary 3.7 above appeared in \cite[Theorem 9(ii)]{5}.\vskip 2mm

In many earlier works, the Bernardi integral operator
$$F(z)^\alpha=\frac{\alpha+c}{z^c}\int_0^zt^{c-1}f(z)^\alpha dt,\;\;\alpha+c>0\eqno(3.1)$$
has received much attention (see \cite{1,2,7,8}). The problem is to study the closure
property of classes of functions under the operator. The converse problem of
determining the largest circle in which $f(z)$ defined by the integral (3.1) belongs to a
certain class of function, given that the function $F(z)$ is a member of the class, has
also been considered by some authors (see \cite{8}). Our concern here is to generalize
and improve those earlier results. Our result is stated as follows.\vskip 2mm

{\bf Theorem 3.9.}\vskip 2mm

{\em Let $F(z)\in T_n^\alpha(\beta)$ and $\alpha+c>0$. Then $f$ given by (3.1) is in $T_n^\alpha(\beta)$ provided $|z|<r_0(\alpha,c)$ where $r_0(\alpha,c)$ is given by

$$r_0(\alpha,c)=\frac{(1+(\alpha+c)^2)^\frac{1}{2}-1}{(\alpha+c)}\eqno(3.2).$$

The result is sharp}.\vskip 2mm

{\bf Proof.}\vskip 2mm

From (3.1), we get
$$(\alpha+c)f(z)^\alpha=DF(z)^\alpha+cF(z)^\alpha.\eqno(3.3)$$

Hence we have
$$(\alpha+c)\frac{D^nf(z)^\alpha}{\alpha^nz^\alpha}=\frac{D^{n+1}F(z)^\alpha}{\alpha^nz^\alpha}+c\frac{D^nF(z)^\alpha}{\alpha^nz^\alpha}.\eqno(3.4)$$

Since $F(z)\in T_n^\alpha(\beta)$, there exists $p\in P$ such that
$$\frac{D^nF(z)^\alpha}{\alpha^nz^\alpha}=\beta+(1-\beta)p(z).\eqno(3.5)$$

Therefore,
$$\frac{D^{n+1}F(z)^\alpha}{\alpha^nz^\alpha}=\alpha(\beta+(1-\beta)p(z))+(1-\beta)zp'(z).\eqno(3.6)$$

Using (3.5) and (3.6) in (3.4), we have
$$(\alpha+c)\left(\frac{D^nf(z)^\alpha}{\alpha^nz^\alpha}-\beta\right)=(1-\beta)[(\alpha+c)p(z)+zp'(z)].\eqno(3.7)$$

Applying Lemma 2.5 on (3.7), we obtain
$$(\alpha+c)\left(Re\frac{D^nf(z)^\alpha}{\alpha^nz^\alpha}-\beta\right)=(1-\beta)\left(\frac{(\alpha+c)(1-r^2)-2r}{1-r^2}\right)Re\;p(z)\eqno(3.8)$$
which implies Re$\tfrac{D^nf(z)^\alpha}{\alpha^nz^\alpha}>\beta$, provided $|z|<r_0(\alpha,c)$.\vskip 2mm

The functions defined by $\tfrac{D^nf(z)^\alpha}{\alpha^nz^\alpha}=\tfrac{1+(1-2\beta)z}{1-z}$ show that the result is sharp. This completes the proof.\vskip 2mm

{\bf Remark 3.10.}\vskip 2mm

We improve Singh's result \cite[Theorem 4']{8} by setting $n=1$ and $\beta=0$. The parameter $\alpha$ and $c$ may not necessarily be integers and $c$ may be as small as `almost' $-\alpha$.

\medskip

{\bf 4.0 Conclusion}\vskip 2mm

In this study, we have provided a sufficient condition for analytic functions normalized by $f(0)=0$ and $f'(z)=1$ to belong to the subclasses $B_n(\alpha)$, which are known to consist of analytic and univalent functions having logarithmic growth in the unit disk \cite{1}. Also we proved another sufficient condition for such analytic function to be of order $\tfrac{1}{2}-B_n(\alpha)$ functions.
For $n=1$, both results are sufficient for univalence in the unit disk. Furthermore, we solved the radius problem associated with the Bernardi integral in relation to a more general class of functions, that is, $T_n^\alpha(\beta)$.\vskip 2mm 

The results of this work generalize and extend many known ones as we have noted in the corollaries and remarks following each of them.


\bigskip

\begin{thebibliography}{99}

\bibitem {1} S. Abdulhalim, \textit{On a class of analytic functions involving the Salagean differential operator}, Tamkang J. Math., {\bf23} (1) (1992), 51-58.

\bibitem {2} K. O. Babalola and T. O. Opoola, \textit{Iterated integral transforms of Caratheodory functions and their applications to analytic and univalent functions}, To appear.

\bibitem {3} K. O. Babalola and T. O. Opoola, \textit{Radius problems for a certain family of analytic functions defined by the Salagean derivative}, Submitted.

\bibitem {4} I. E. Bazilevic, \textit{On a case of integrability by quadratures of the equation of Loewner-Kufarev}, Mat. Sb. {\bf37} (79), (1955), 471-476 (Russian).

\bibitem {5} S. S. Miller and P. T. Mocanu, \textit{Second order differential inequalities in the complex plane}, J. Math. Anal. Appl. {\bf65}, (1978), 289-305.

\bibitem {6} J. W. Noonan, \textit{On close-to-convex functions of order $\beta$}, Pacific J. Math. {\bf44} (1), (1973), 263-280.

\bibitem {7} T. O. Opoola, \textit{On a new subclass of univalent functions}, Mathematica (Cluj) {\bf36}, 59 No. 2 (1994), 195-200.

\bibitem {8} R. Singh, \textit{On Bazilevic functions}, Proc. Amer. Math. Soc. {\bf38} (1973), 261-271.

\bibitem {9} K. Yamaguchi, \textit{On functions satisfying $Re\{f(z)/z\}>0$}, Proc. Amer. Math. Soc. {\bf17} (1966),
588-591.

\end{thebibliography}
\end{document}